\documentclass[a4paper,12pt]{amsart}
\usepackage{amssymb,amsmath,a4wide,graphicx,stmaryrd,fullpage,setspace,microtype,textcomp,tikz,tikz-cd,accents,mathtools,abraces}
\usepackage{hyperref}
\usepackage[shortlabels]{enumitem}
\definecolor{dark-red}{rgb}{0.5,0.15,0.15}
\hypersetup{
	colorlinks   = true, %Colours links instead of ugly boxes
	urlcolor     = dark-red, %Colour for external hyperlinks
	linkcolor    = black, %Colour of internal links
	citecolor   = dark-red %Colour of citations
}
\usepackage[T1]{fontenc}
\usepackage{lmodern}
\usepackage[numbers,sort&compress]{natbib}

\title{Branching space of precubical set}

\author[P. Gaucher]{Philippe Gaucher}

\address{Universit\'e Paris Cit\'e, CNRS, IRIF, F-75013, Paris, France}

\urladdr{http://www.irif.fr/{\~{}}gaucher} 

\makeatletter
\@namedef{subjclassname@2020}{%
	\textup{2020} Mathematics Subject Classification}
\makeatother


\swapnumbers

\makeatletter
\let\P\@undefined
\let\leq\@undefined
\let\geq\@undefined
\let\top\@undefined
\makeatother
\newcommand{\K}{\mathcal{K}}
\newcommand{\p}{\times}
\renewcommand{\vec}{\overrightarrow}
\newcommand{\P}{\mathbb{P}}
\DeclareMathOperator{\hop}{ho\P}

\newcommand{\dtop}{{\brm{Flow}}}
\newcommand{\set}{{\brm{Set}}}
\newcommand{\glob}{{\rm{Glob}}}

\newcommand{\top}{{\mathbf{Top}}}
\newcommand{\iso}{\cong}
\newcommand{\vI}{\vec{I}}
\newcommand{\leq}{\leqslant}
\newcommand{\geq}{\geqslant}
\newcommand{\tr}[1]{{\langle{#1}\rangle}}

\newcommand{\brm}[1]{\rm{\mathbf{#1}}}
\newcommand{\ttop}{{\brm{TOP}}}

\DeclareMathOperator{\id}{Id}
\newcommand{\liminj}{\varinjlim}

\newcommand{\ptop}[1]{{\brm{{#1}dTop}}}
\newcommand{\cocartesian}{\arrow[lu, phantom, "\ulcorner"{font=\Large}, pos=0]}

\DeclareMathOperator{\im}{im}

\newcommand{\de}{\partial}
\newcommand{\rest}{\!\!\restriction\!\!}
\DeclareMathOperator{\Sub}{Sub}

\newtheorem*{thmN}{Theorem}

\newtheorem{thm}{Theorem}[section]
\newtheorem{prop}[thm]{Proposition}
\newtheorem{lem}[thm]{Lemma}
\newtheorem{cor}[thm]{Corollary}
\newcommand{\bth}{\begin{thm}}
\renewcommand{\eth}{\end{thm}}
\newcommand{\bpf}{\begin{proof}}
\newcommand{\epf}{\end{proof}}
\newcommand{\bc}{\begin{cor}}
\newcommand{\ec}{\end{cor}}
\theoremstyle{definition}
\newtheorem{defn}[thm]{Definition}

\newcommand{\bd}{\begin{defn}}
\newcommand{\ed}{\end{defn}}
\newtheorem{nota}[thm]{Notation}

\makeatletter
\def\varholim@#1#2{%
	\vtop{\m@th\ialign{##\cr
			\hfil$#1\operator@font holim$\hfil\cr
			\noalign{\nointerlineskip\kern1.5\ex@}#2\cr
			\noalign{\nointerlineskip\kern-\ex@}\cr}}%
}
\def\holimproj{%
	\mathop{\mathpalette\varholim@{\leftarrowfill@\textstyle}}\nmlimits@
}
\def\holiminj{%
	\mathop{\mathpalette\varholim@{\rightarrowfill@\textstyle}}\nmlimits@
}
\makeatother

\makeatletter
\newcommand*{\@opargbegintheorem}[3]{\trivlist
	\item[\hskip \labelsep{\bfseries #1\ #2}] \textbf{(#3)}\ \itshape}
\makeatother

\allowdisplaybreaks

\begin{document}

\begin{abstract} 
	Using the notion of short natural directed path, we introduce the homotopy branching space of a precubical set. It is unique only up to homotopy equivalence. We prove that, for any precubical set, it is homotopy equivalent to the branching space of any q-realization, any m-realization and any h-realization of the precubical set as a flow. As an application, we deduce the invariance of the homotopy branching space and of the branching homology up to cubical subdivision. By reversing the time direction, the same results are obtained for the merging space and the merging homology of a precubical set.
\end{abstract}

\maketitle
\tableofcontents
\hypersetup{linkcolor = dark-red}

\section{Introduction}

\subsection*{Presentation}

The primary goal of this paper is to present a cubical analogue of the results of \cite{topological-branching}, which treats the globular setting. Specifically we introduce the notion of branching for precubical sets and proving its invariance by cubical subdivisions. 

Directed Algebraic Topology (DAT) aims to model concurrent systems using tools from algebraic topology \cite{DAT_book}. The core idea is that an $n$-cube represents the concurrent execution of $n$ actions. The geometric realization $[0,1]^n$ is equipped with the set of of all continuous paths which are nondecreasing with respect to each axis of coordinates. They represent all possible executions. The nondecreasing condition reflects the irreversibility of time. In this framework, precubical sets (Definition~\ref{def:precubical-set}) generalize the notion of cube, allowing them to assemble into complex structures that capture the behavior of concurrent systems. While precubical sets are general enough to encompass all examples from concurrency theory, most computationally relevant precubical sets are non-positively curved in the sense of \cite[Definition~1.28 and Proposition~1.29]{zbMATH07226006} or, at worst, proper as defined in \cite[page~499]{MR3722069}.

The notions of branching space and branching homology (along with their duals, merging space and merging homology) are relevant in DAT because they encode local information about the causal structure of a concurrent system. They detect the non-deterministic branching and merging areas of execution paths. Originally, these concepts were introduced not for precubical sets but for flows (Definition~\ref{def:flow}) \cite{exbranching}: see Definitions~\ref{def:branching-flow} and \ref{def:hombrdef-flow}.

A realization functor from precubical sets to flows was introduced in \cite[Definition~7.2]{ccsprecub}, defined as the unique colimit-preserving functor mapping the $n$-cube to $(\{0<1\}^n)^{cof}$ where $\{0<1\}^n$ is the product poset of $n$ copies of $\{0<1\}$, viewed as a flow, and where $(-)^{cof}$ is a fixed q-cofibrant functor of the q-model structure of flows (see Theorem~\ref{thm:three}). This functor is an example of a q-realization functor in the sense of Definition~\ref{def-rea-flow} and is unique only up to homotopy by \cite[Theorem~3.8]{NaturalRealization}, as it depends on the arbitrary choice of a q-cofibrant replacement functor in the category of flows.

The purpose of this paper is to bypass this indirect construction. We provide a direct definition of the branching space and branching homology for precubical sets without relying on any realization functor and prove that it recovers the previous definition.

When this paper was initially conceived, our intention was to adapt the notion of germ of short directed paths, originally introduced in the globular setting \cite{topological-branching}, to the cubical setting. However, the cubical setting introduced unexpected challenges:
\begin{itemize}[leftmargin=*]
	\item Fixed Natural Length: Germs of short directed paths are replaced by natural directed paths with a \textit{fixed} natural length $\epsilon\in ]0,1[$.
	\item Inadequacy of Germs: The notion of germ of directed paths is unsuited to the cubical setting by Proposition~\ref{prop:no-germs}. Instead, Proposition~\ref{prop:ho-germ} shows that elements of the branching space should be interpreted as homotopy germs of short natural directed paths.
	\item Uniqueness Up to Homotopy: Unlike the globular case, where the branching space is unique up to homeomorphism, the branching space of a precubical set is unique only up to homotopy equivalence. This difference arises because the realization of a precubical set as a directed space (Definition~\ref{rea}) contains more directed paths than its realization in the globular setting of \cite{topological-branching}. This geometric richness is the underlying reason for Proposition~\ref{prop:no-germs}. 
\end{itemize}

We now state the main results of the paper. Let $\mathcal{P}^-(K,\epsilon)$ be the space of natural directed paths of natural length $\epsilon\in ]0,1[$ of a precubical set $K$ starting from a vertex of $K$ (Definition~\ref{def:branching-space}).

\begin{thmN} (Theorem~\ref{thm:G-colimit-preserving})
	The functor $\mathcal{P}^-(-,\epsilon):\square^{op}\set \to \top$ from precubical sets to topological spaces is colimit-preserving.
\end{thmN}

\begin{thmN} (Theorem~\ref{thm:P-homotopy-invariant})
	Let $\epsilon,\epsilon'\in ]0,1[$. For all precubical sets $K$, there is a natural homotopy equivalence between m-cofibrant spaces 
	\[
	\mathcal{P}^-(K,\epsilon) \simeq \mathcal{P}^-(K,\epsilon').
	\]
\end{thmN}

\begin{thmN} (Corollary~\ref{cor:top-flow-homotopy-equivalent})
	Let $r\in \{q,m,h\}$. Let $F:\square^{op}\set\to \dtop$ be an r-realization functor. Then for all precubical sets $K$, there is a natural homotopy equivalence 
	\[
	\mathcal{P}^-(K,\epsilon) \simeq \P^-F(K)
	\]
	where $\P^-:\dtop \to \top$ is the branching space functor from flows to topological spaces defined in Theorem~\ref{thm:def-Cminus}.
\end{thmN}

Finally, we obtain as an application the invariance by cubical subdivisions:

\begin{thmN} (Theorem~\ref{thm:sub})
	Let $K,L$ be two precubical sets. Consider a cubical subdivision $f:K\dashrightarrow L$ ($K_0$ is therefore identified with a subset of $L_0$) in the sense of Definition~\ref{def:cubical-sbd}. Let $||-||:\square^{op}\set\to \dtop$ be an r-realization functor with $r\in \{q,m,h\}$. For all $\alpha\in K_0$. there is a homotopy equivalence $\P^-_\alpha ||K|| \simeq \P^-_\alpha ||L||$ where $\P^-_\alpha X$ is the branching space of the flow $X$ at $\alpha$. For all $\alpha\in L_0\backslash K_0$, the space $\P^-_\alpha||L||$ is contractible: cubical subdivisions do not create non-deterministic branching areas of execution paths.
\end{thmN}

\begin{thmN} (Corollary~\ref{cor:enfinfin})
	Let $K,L$ be two precubical sets. Let $f:K\dashrightarrow L$ be a cubical subdivision. Then for any $n\geq 0$, there is an isomorphism $H_n^-(||K||)\iso H_n^-(||L||)$ for any r-realization functor $||-||:\square^{op}\set\to \dtop$.
\end{thmN}

A generating subdivision in the category of flows is defined as the pushout of a q-cofibrant replacement of an inclusion of finite bounded posets $P_1\subset P_2$, viewed as a map of flows, that preserves the smallest and the biggest elements \cite[Definition~8.3]{MultipointedSubdivision}. If a cubical analogue of the main results of \cite{MultipointedSubdivision} were established, we could prove that, for any cubical subdivision $f:K\dashrightarrow L$, the induced map of flows $||f||:||K||\to ||L||$ is a zigzag of transfinite compositions of generating subdivisions and of weak equivalences of flows. Then Corollary~\ref{cor:enfinfin} would follow directly from \cite[Corollary~11.3]{3eme}. Indeed, the latter corollary guarantees the invariance of branching and merging homology theories under transfinite compositions of generating subdivisions. We intend to explore such a cubical analogue in future work.

A key takeaway from these results and those in \cite{topological-branching} is that the concept of a natural directed path seems essential for defining a meaningful branching space. While the term ``natural length'' is not explicitly used in \cite{topological-branching}, it is implicitly present: the short directed paths, as defined in \cite[Definition~4.2]{topological-branching}, correspond precisely to directed paths with globular natural lengths in the open interval $]0,1[$.

This suggests a deeper issue: it is likely impossible to construct a coherent notion of branching space within the category of directed spaces alone. The category is simply too broad and lacks the necessary structure. One potential approach would be to enrich directed spaces with a notion of natural length. 

For every precubical set $K$, its directed realization $\vec{|K|}$ as a directed space (Definition~\ref{rea}) can be equipped by \cite[Section~2]{MR2521708} with a map $p:\vec{|K|}\to \vec{\mathbf{S}}^1$ where $\vec{\mathbf{S}}^1$ is the directed circle with $p(|c|_{geom}(x_1,\dots,x_n))=x_1+\dots + x_n \mod 1$ for an $n$-cube $c$ of $K$. A possible solution could therefore involve working in a comma category, where a short directed path of $X\to \vec{\mathbf{S}}^1$ would be defined as a directed path of the directed space $X$ such that the image by the map $X\to \vec{\mathbf{S}}^1$ does not cover the entire circle. However, this remains highly speculative. Indeed, the globular setting does not fit into this formal setting. By \cite[Corollary~4]{Moore3}, unlike the cubical natural length, the globular natural length is not constant by continuous deformations of directed paths preserving the extremities. Instead, it takes finitely many values along.

As a concluding observation, the constructions developed in this paper can be transposed to the setting of merging space and merging homology by inverting the temporal orientation. The reader may verify the details along the same lines.

\subsection*{Outline of the paper}

Section~\ref{section:reminders} collects some reminders about general topology and homotopy theory. It is also useful to fix the notations and provides some references for the reader.
	
Section~\ref{section:defs} recalls the notion of directed path in the geometric realization of a precubical set, the notion of directed space, and the realization of a precubical set as a directed space.
	
Section~\ref{section:homotopy-branching} introduces the notion of short natural directed path. And then it introduces the homotopy branching space of a precubical set. Some basic topological properties related to the $n$-cube and its boundary are proved.
	
Section~\ref{section:cat-branching} explores the categorical properties of the homotopy branching space functor. It culminates with Theorem~\ref{thm:G-colimit-preserving} which establishes that the homotopy branching space functor is colimit-preserving.
	
Section~\ref{section:homotopic-branching} explores the homotopical properties of the homotopy branching space functor. It concludes with the proof in Theorem~\ref{thm:P-homotopy-invariant} that the choice of the parameter $\epsilon$ is not important provided that we work up to homotopy equivalence. Proposition~\ref{prop:ho-germ} explains why the homotopy branching space can be viewed as the space of homotopy germs of short natural paths.
	
Section~\ref{section:flow-branching} recalls the notion of homotopy branching space of a flow. Then it is proved in Theorem~\ref{thm:top-flow-homotopy-equivalent} that the two notions, the topological one of Section~\ref{section:homotopy-branching} and the one of this section, coincide up to homotopy equivalence. Finally Corollary~\ref{cor:top-flow-homotopy-equivalent} establishes that the q-realization functor can be replaced by any r-realization functor for $r\in \{q,m,h\}$.
	
As an application, it is proved in Section~\ref{section:sub} that the homotopy branching space and the branching homology are invariant up to cubical subdivision in Theorem~\ref{thm:sub} and Corollary~\ref{cor:enfinfin}. It is the cubical analogue of \cite[Theorem~7.8]{topological-branching} and \cite[Corollary~7.9]{topological-branching}.

\section{Topological and homotopical reminders}
\label{section:reminders}

In this paper, compact means Hausdorff quasi-compact. The category $\top$ denotes either the category of \textit{$\Delta$-generated spaces} or of \textit{$\Delta$-Hausdorff $\Delta$-generated spaces} (cf. \cite[Section~2 and Appendix~B]{leftproperflow}). It is Cartesian closed by a result due to Dugger and Vogt recalled in \cite[Proposition~2.5]{mdtop} and locally presentable by \cite[Corollary~3.7]{FR}. The internal hom is denoted by $\ttop(-,-)$. The right adjoint $k_\Delta:\mathcal{TOP}\to \top_\Delta$ of the inclusion from $\Delta$-generated spaces to general topological spaces is called the $\Delta$-kelleyfication. Let $\phi_U:k_\Delta(U)\to U$ be the counit. It preserves the underlying sets. One has $\ttop(-,-)=k_\Delta(\ttop_{co}(-,-))$ where $\ttop_{co}(-,-)$ means the set of continuous maps equipped with the compact-open topology. Every open subset of a $\Delta$-generated space equipped with the relative topology is $\Delta$-generated. A quotient map is a continuous map $f:X\to Y$ of $\Delta$-generated spaces which is onto and such that $Y$ is equipped with the final topology. The space $Y$ is called a final quotient of $X$. Every $\Delta$-generated space is sequential. The category $\top$ is equipped with its q-model structure. The m-model structure \cite{mixed-cole} and the h-model structure \cite{Barthel-Riel} of $\top$ are also used in various places of the paper. The q-model structure of $\top$ is combinatorial. The m-model structure and the h-model structure of $\top$ are accessible (e.g. see \cite[Proposition~4.2]{QHMmodel}).

A \textit{closed inclusion} $f:A\to B$ of $\Delta$-generated spaces is a one-to-one continuous map such that $f(A)$ is a closed subset of $B$ and such that $f$ induces a homeomorphism between $A$ and $f(A)$ equipped with the relative topology. A one-to-one map of $\Delta$-generated spaces $f:A\to B$ is a \textit{$\Delta$-inclusion} if for all $\Delta$-generated spaces $Z$, the set map $Z\to A$ is continuous if and only if the composite set map $Z\to A\to B$ is continuous. Every closed inclusion of $\Delta$-generated spaces is a $\Delta$-inclusion with closed image. The converse is false: see \cite[Section~2]{topological-branching}.

\bth \label{thm:recognizing-closed-inclusions-1} \cite[Theorem~2.4]{topological-branching}
Let $f:A\to B$ be a one-to-one continuous map between $\Delta$-generated spaces. Suppose that $B$ is $\Delta$-Hausdorff. The map $f$ is a closed inclusion if and only if for any sequence $(x_n)_{n\geq 0}$ of $A$, if the sequence $(f(x_n))_{n\geq 0}$ of $B$ is convergent, so is the sequence $(x_n)_{n\geq 0}$ of $A$.
\eth

Let $[n] = \{0,1\}^n$ for $n \geq 1$. Let $\{0,1\}^0=[0,1]^0=[0]=\{()\}$. Let $0_n=(0,\dots,0)\in [n]$ and $1_n=(1,\dots,1)\in [n]$. Let $\delta_i^\alpha : [0,1]^{n-1} \rightarrow [0,1]^n$ be the continuous map defined for $1\leq i\leq n$ and $\alpha \in \{0,1\}$ by $\delta_i^\alpha(\epsilon_1, \dots, \epsilon_{n-1}) = (\epsilon_1,\dots, \epsilon_{i-1}, \alpha, \epsilon_i, \dots, \epsilon_{n-1})$. The small category $\square$ is the subcategory of the category of sets with the set of objects $\{[n],n\geq 0\}$ and generated by the \textit{coface maps} $\delta_i^\alpha$. It is called the \textit{box category}.

\bd \label{def:precubical-set} A \textit{precubical set} $K$ is a presheaf over $\square$. The category of precubical sets is denoted by $\square^{op}\set$. 
\ed

An element $c$ of $K_n$ is called an \textit{$n$-cube} and we set $n=\dim(c)$. For all $n\geq 1$, there are $2n$ \textit{face maps} $\de_i^\alpha:K_n\to K_{n-1}$. the $n$-cube is the precubical set $\square[n]=\square(-,[n])$ for $n\geq 0$. For a precubical set $K$, the precubical set $K_{\leq p}$ is the precubical subset of $K$ containing only the cube of $K$ of dimension lower than $p$. The boundary of $\square[n]$ is the precubical set $\square[n]_{\leq n-1}$. It is denoted by $\de \square[n]$. In particular, one has $\de \square[0] = \varnothing$. 

The \textit{initial vertex} of an $n$-cube $c$ of a precubical set $K$ is the vertex $c^-=\de_1^0\dots\de_n^0c$. When $c\in K_0$, one has $c^-=c$. For $\alpha\in K_0$, let $\mathcal{C}^-_\alpha(K)=\{c\in K\mid \dim(c)\geq 1 \hbox{ and }c^-=\alpha\}$.

A cocubical object $F:\square\to\K$ of a cocomplete category $\K$ gives rise to a colimit-preserving functor 
\[
\widehat{F}(K)=\int^{[n]\in \square} K_n.F(\square[n]).
\]
By e.g. \cite[Remark~3.2.3]{coend-calculus} \cite[Proposition~2.3.2]{realization}, this mapping induces an equivalence of categories between the category of cocubical objects of $\K$ and the colimit-preserving functor from $\square^{op}\set$ to $\K$. 

The box category $\square$ is a direct Reedy category with the degree function $d([n])=n$. Let $\mathcal{M}$ be a model category. Equip the category of cocubical objects $\mathcal{M}^\square$ with its Reedy model structure. It coincides with the projective model structure since the Reedy category $\square$ is direct. 

\bth \label{thm:tout}
Let $\mathcal{M}$ be a simplicial model category. 
\begin{enumerate}
	\item A cocubical object $F:\square\to \mathcal{M}$ is Reedy cofibrant if and only if for all $n\geq 0$, the map $\widehat{X}(\de\square[n]\subset \square[n])$ is a cofibration of $\mathcal{M}$.
	\item For any Reedy cofibrant cocubical object $F:\square\to \mathcal{M}$, for all precubical sets $K$, $\widehat{F}(K)$ is cofibrant.
	\item Consider two Reedy cofibrant cocubical objects $F,G:\square\to \mathcal{M}$ such that there exists a cocubical object $I:\square\to \mathcal{M}$ and two trivial Reedy fibrations (i.e. trivial objectwise fibrations) \[
	\begin{tikzcd}[row sep=small,column sep=2em]
		F \arrow[r,twoheadrightarrow,"\simeq"] & I & G \arrow[l,twoheadrightarrow,"\simeq"']
	\end{tikzcd}
	\]
	then there exists a natural transformation $\mu:\widehat{F}\Rightarrow \widehat{G}$ such that for all precubical sets $K$, there is a homotopy equivalence $\mu_K:\widehat{F}(K) \simeq \widehat{G}(K)$.
\end{enumerate}  
\eth

\bpf
(1) is \cite[Proposition~2.3.1]{realization}. (2) is an easy induction. (3) is \cite[Theorem~2.3.3]{realization}.
\epf

A subset $A$ of a topological space $X$ is called a \textit{retract} of $X$ if there is a map $r:X\to A$ such that $r\rest A=\id$, which is called a retraction. It is well-known that every retract of a space $X$ is closed in $X$ whenever the diagonal of $X$ is closed in $X\p X$. It is the case e.g. when $X$ is metrizable or when $X$ is $\Delta$-Hausdorff. A \textit{neighborhood retract} of $X$ is a closed subset of $X$ that is a retract of some neighborhood in $X$. A metrizable space $X$ is called an \textit{absolute neighborhood retract} (ANR) (resp. an \textit{absolute retract} (AR)) if $X$ is a neighborhood retract (or a retract) of an arbitrary metrizable space that contains $X$ as a closed subspace. A modern reference for AR and ANR is \cite[Chapter~6]{Sakai}. These notions are used in the proof of Proposition~\ref{prop:cubical-cofibration-1}.

\section{From precubical set to directed space}
\label{section:defs}

Let $U$ be a topological space. A \textit{(Moore) path} of $U$ consists of a continuous map $\gamma:[0,\ell]\to U$. The real number $\ell$ is called the \textit{length} of $\gamma$. Let $\gamma_1:[0,\ell_1]\to U$ and $\gamma_2:[0,\ell_2]\to U$ be two Moore paths of a topological space $U$ such that $\gamma_1(\ell_1)=\gamma_2(0)$. The \textit{Moore composition} $\gamma_1*\gamma_2:[0,\ell_1+\ell_2]\to U$ is the Moore path defined by 
\[
(\gamma_1*\gamma_2)(t)=
\begin{cases}
	\gamma_1(t) & \hbox{ for } t\in [0,\ell_1]\\
	\gamma_2(t-\ell_1) &\hbox{ for }t\in [\ell_1,\ell_1+\ell_2].
\end{cases}
\]
The Moore composition of Moore paths is strictly associative. Let $\gamma_1$ and $\gamma_2$ be two continuous maps from $[0,1]$ to some topological space such that $\gamma_1(1)=\gamma_2(0)$. The continuous map defined by 
\[
(\gamma_1 *_N \gamma_2)(t) = \begin{cases}
	\gamma_1(2t)& \hbox{ if }0\leq t\leq \frac{1}{2},\\
	\gamma_2(2t-1)& \hbox{ if }\frac{1}{2}\leq t\leq 1
\end{cases}
\]
is called the \textit{normalized composition}.

The mapping $[n]\mapsto [0,1]^n$ yields a well-defined \textit{cocubical object} of $\top$, i.e. a functor from $\square$ to $\top$. The latter gives rise to a colimit-preserving functor from precubical sets to topological spaces called the \textit{geometric realization}
\[
|K|_{geom} = \int^{[n]\in \square} K_n.[0,1]^n.
\]
The topological space $|K|_{geom}$ is always a CW-complex. In particular it is always Hausdorff (see e.g. \cite[Proposition~3.12]{DirectedDegeneracy}).

\begin{nota}
	The set $\mathcal{I}$ of non-decreasing continuous maps from $[0,1]$ to $[0,1]$ is equipped with the $\Delta$-kelleyfication of the relative topology induced by the inclusion $\mathcal{I} \subset \ttop([0,1],[0,1])$.
\end{nota}

\bd \cite[Definition~1.1]{mg} \cite[Definition~4.1]{DAT_book} \label{def:directed_space}
A \textit{directed space} is a pair $X=(|X|,\vec{P}(X))$ consisting of a topological space $|X|$ and a set $\vec{P}(X)$ of continuous paths from $[0,1]$ to $|X|$ satisfying the following axioms:
\begin{itemize}
	\item $\vec{P}(X)$ contains all constant paths;
	\item $\vec{P}(X)$ is closed under normalized composition;
	\item $\vec{P}(X)$ is closed under reparametrization by an element of $\mathcal{I}$.
\end{itemize}
The space $|X|$ is called the \textit{underlying topological space} or the \textit{state space}. The elements of $\vec{P}(X)$ are called \textit{directed paths}. A morphism of directed spaces is a continuous map between the underlying topological spaces which takes a directed path of the source to a directed path of the target. The category of directed spaces is denoted by $\ptop{}$. The space $\vec{P}(X)$ is equipped with the $\Delta$-kelleyfication of the compact-open topology.
\ed

The category $\ptop{}$ is locally presentable and in particular cocomplete (see \cite[Theorem~4.2]{FR} and the remark before \cite[Proposition~3.6]{GlobularNaturalSystem}). By \cite[Remark~4.3~(ii)]{FR}, the set of directed paths of a colimit of directed spaces is given by a final structure: it is obtained by taking the closure by concatenation and reparametrization. Thus, the functor from $\ptop{}$ to $\top$ taking a directed space $X=(|X|,\vec{P}(X))$ to its state space $|X|$ is topological in the sense of \cite[Chapter~VI]{topologicalcat}. In particular, it is colimit-preserving.

The geometric realization $|\square[n]|_{geom}=[0,1]^n$ can be equipped with a structure of directed space as follows. The directed paths are the continuous paths $\gamma=(\gamma_1,\dots,\gamma_n):[0,1]\to [0,1]^n$ such that each $\gamma_i$ belongs to $\mathcal{I}$ for $i=1,\dots,n$. We obtain a directed space denoted by $\vec{[0,1]^n}$. The mapping $[n]\mapsto \vec{[0,1]^n}$ gives rise to a cocubical object of $\ptop{}$, the coface maps preserving the local ordering. We obtain a colimit-preserving functor from precubical sets to directed spaces:

\bd \label{rea}
Let $K$ be a precubical set. The \textit{directed realization} $\vec{|K|}$ of $K$ is the directed space \[\vec{|K|}=\int^{[n]\in \square} K_n.\vec{[0,1]^n}.\]
\ed 

The state space of the directed space $\vec{|K|}$ is $|K|_{geom}$. A directed path of $[0,1]^n$ of length $\ell$ is a continuous map $\gamma:[0,\ell] \to [0,1]^n$ which is non-decreasing with respect to each axis of coordinates. A directed path of an $n$-cube $c\in K_n$ of length $\ell \geq 0$ is a composite continuous map of the form $|c|_{geom}\gamma:[0,\ell] \to [0,1]^n\to |K|_{geom}$ such that $\gamma:[0,\ell]\to [0,1]^n$ is a directed path of length $\ell\geq 0$ ($c$ giving rise to a map $\square[n]\to K$ using the Yoneda lemma). A directed path in $|K|_{geom}$ of length $\ell\geq 0$ is a continuous path $[0,\ell] \to |K|_{geom}$ which is a Moore composition of the form $(|c_1|_{geom}\gamma_1) * \dots *(|c_p|_{geom}\gamma_p)$ with $p\geq 1$ such that $\ell=\ell_1+\dots+\ell_p$ where $\ell_i$ is the length of $\gamma_i$ for $1\leq i \leq p$. The choice of $c_1,\dots,c_p$ is not unique.

\section{Homotopy branching space of a precubical set}
\label{section:homotopy-branching}

A directed path $\gamma=(\gamma_1,\dots,\gamma_n)$ from $[0,\ell]$ to $[0,1]^n$ starting from $0_n$ is \textit{natural} if $t=\gamma_1(t)+\dots + \gamma_n(t)$ for all $t\in [0,\ell]$ (see \cite[Section~2.2]{MR2521708} and \cite[Definition~4.8]{NaturalRealization}). There is a more general notion of natural directed path introduced in \cite[Section~2.2]{MR2521708} for the directed space $\vec{|K|}$ associated with a precubical set $K$. In the language of Lawvere metric spaces \cite{LawvereMetric}, a natural directed path is a Moore directed path which is a Moore composition of quasi-isometry. This is the point of view adopted in \cite{DirectedDegeneracy}. We do not need in this paper the notion of natural directed path in full generality.  

\bd 
Let $n\geq 1$. A \textit{short natural directed path} of $[0,1]^n$ is a natural directed path $\phi:[0,\epsilon]\to [0,1]^n$ such that $\phi(0)=0_n$ and $0<\epsilon<1$. The real number $\epsilon$ is called the \textit{natural length} of the short natural directed path. The set of short natural directed paths of $[0,1]^n$ of length $\epsilon$ is denoted by $N_n(\epsilon)$. 
\ed

\bd \label{def:branching-space}
Let $\epsilon\in ]0,1[$. Let $K$ be a precubical set. The \textit{homotopy branching space} $\mathcal{P}^-(K,\epsilon)$ of $K$ is the space 
\[
\mathcal{P}^-(K,\epsilon) = \coprod_{\alpha\in K^0} \mathcal{P}^-_\alpha(K,\epsilon).
\]
such that 
$
\mathcal{P}^-_\alpha(K,\epsilon) = \big\{|c|_{geom}\phi\mid c\in \mathcal{C}^-_\alpha(K)\hbox{ and }\phi\in N_{\dim(c)}(\epsilon)\big\}
$
is equipped with the $\Delta$-Kelleyfication of the relative topology with respect to the one of  $\ttop([0,\epsilon],|K|_{geom})$. 
\ed

The reason of the terminology is that the space $\mathcal{P}^-(K,\epsilon)$ is unique up to homotopy by Theorem~\ref{thm:P-homotopy-invariant}. So the choice of $\epsilon\in ]0,1[$ does not matter. Further comments about Definition~\ref{def:branching-space} are given in Proposition~\ref{prop:no-germs} and Proposition~\ref{prop:ho-germ}.

\begin{prop}
	The mapping $K\mapsto \mathcal{P}^-(K,\epsilon)$ induces a functor from precubical sets to topological spaces.
\end{prop}

\bpf Let $f:K\to L$ be a map of precubical sets. For all $c\in\mathcal{C}^-_\alpha(K)$, $f(c)\in \mathcal{C}^-_{f(\alpha)}(L)$ and $\dim(c)=\dim(f(c))$. Hence the proof is complete since $0<\epsilon<1$. 
\epf

Proposition~\ref{prop:no-germs} proves the unsuitability of the notion of germ in the cubical setting.

\begin{prop} \label{prop:no-germs} (compare with Proposition~\ref{prop:cubical-cofibration-1})
	Let $K$ be a precubical set. Let $\alpha\in K_0$. Introduce the quotient $\mathcal{G}^-_\alpha(K)$ of the homotopy branching space $\mathcal{P}^-_\alpha(K,\epsilon)$ by the equivalence relation $\gamma_1\simeq^- \gamma_2$ if there exists $0<\epsilon'<\epsilon$ such that $\gamma_1\rest_{[0,\epsilon']} = \gamma_2\rest_{[0,\epsilon']}$. Then the map $\mathcal{G}^-_{0_2}(\de\square[2])\to \mathcal{G}^-_{0_2}(\square[2])$ is not a closed inclusion, and therefore not an r-cofibration for $r\in \{q,m,h\}$.
\end{prop}

\bpf
Let $\tr{-}^-:\mathcal{P}^-_\alpha(K,\epsilon)\to \mathcal{G}^-_\alpha(K)$ be the quotient map. Chose an integer $N$ such that \[N>\frac{1}{\min(\epsilon,2-\epsilon)}.\] Then $1/N<\min (\epsilon,2-\epsilon)$. We deduce 
\[
\bigg\{\frac{\epsilon-\frac{1}{N}}{2},\frac{\epsilon+\frac{1}{N}}{2}\bigg\} \in ]0,1[.
\] 
Consider the sequence of natural directed paths $(\gamma_{1/n})_{n\geq N}$ depicted in Figure~\ref{fig:not-closed}. For every $n\geq N$, $\tr{\gamma_{1/n}}^-\in \mathcal{G}^-_{0_2}(\de\square[2])$. However, the limit in $\mathcal{G}^-_{0_2}(\square[2])$ does not belong to $\mathcal{G}^-_{0_2}(\de\square[2])$.
\epf 

\begin{figure}
	\def\n{5}
	\begin{tikzpicture}[black,scale=4,pn/.style={circle,inner sep=0pt,minimum width=4pt,fill=dark-red}]
		\fill [color=gray!15] (0,0) -- (0,1) -- (1,1) -- (1,0) -- cycle;
		\draw[dark-red][->][thick] (0,0) -- (0,1/3) -- (1/4,1/4+1/3) ;
		\draw (0,0) node[pn] {} node[black,below left] {$(0\,,0)$};
		\draw (0,1/3) node[pn] {} node[black,left] {$(0\,,h)$};
		\draw (1/4,1/4+1/3) node[pn] {} node[black,right] {$(\frac{\epsilon-h}{2}\,,\frac{\epsilon+h}{2})$};
		\draw (1,1) node[black,above right] {$(1\,,1)$};
		\draw (0,1) node[black,above left] {$(0\,,1)$};
		\draw (1,0) node[black,below right] {$(1\,,0)$};
	\end{tikzpicture}
	\caption{The natural directed path $\gamma_h$ from $(0,0)$ to $(0,h)$ and from $(0,h)$ to $(\frac{\epsilon-h}{2}\,\frac{\epsilon+h}{2})$ of natural length $\epsilon$ in the $2$-cube $[0,1]\p [0,1]$}
	\label{fig:not-closed}
\end{figure}

\begin{lem} (Cutler's observation) \cite[Lemma~5.2]{topological-branching} \label{lem:open-delta-k}
	Let $M$ be a general topological space. If $U$ is an open subset of $M$, then the topology of $k_\Delta(U)$ is the relative topology with respect to the inclusion into $k_\Delta(M)$.
\end{lem}

\begin{prop} \label{prop:P-pushout}
	Consider a pushout diagram of precubical sets 
	\[
	\begin{tikzcd}[row sep=3em,column sep=3em]
		\de\square[n] \arrow[r,"g"] \arrow[d,"\subset"'] & K \arrow[d] \\
		\square[n] \arrow[r,"\widehat{g}"] & \cocartesian L
	\end{tikzcd}
	\]
	For $n=0$ or for $n\geq 1$ and $\alpha\neq g(0_n)$, the map of precubical sets $K\to L$ induces an equality $\mathcal{P}^-_\alpha(K,\epsilon) = \mathcal{P}^-_\alpha(L,\epsilon)$. For $n\geq 1$ and $\alpha=g(0_n)$, the map of precubical sets $K\to L$ induces a closed inclusion $\mathcal{P}^-_\alpha(K,\epsilon) \to \mathcal{P}^-_\alpha(L,\epsilon)$.
\end{prop}

\bpf
When $n=0$, one has $L=K\sqcup \square[0]$. This implies that $\mathcal{C}^-_\alpha(L)=\mathcal{C}^-_\alpha(K)$ for all $\alpha\in K^0$. We obtain the equality $\mathcal{P}^-_\alpha(K,\epsilon) = \mathcal{P}^-_\alpha(L,\epsilon)$. Assume that $n\geq 1$. In this case, the map of precubical sets $K\to L$ induces a bijection of sets $K_0\iso L_0$ between the vertices of $K$ and $L$. If $\alpha\neq g(0_n)$, then $\mathcal{C}^-_\alpha(K)=\mathcal{C}^-_\alpha(L)$. We obtain the equality $\mathcal{P}^-_{\alpha}(K,\epsilon)=\mathcal{P}^-_{\alpha}(L,\epsilon)$ in this case. Assume now that $\alpha= g(0_n)$. In this case, one has $\mathcal{C}^-_\alpha(L)=\mathcal{C}^-_\alpha(K) \sqcup \{\widehat{g}\}$. To prove that the map $\mathcal{P}^-_\alpha(K,\epsilon) \to \mathcal{P}^-_\alpha(L,\epsilon)$ is a closed inclusion, we use the same technique as the one used for the proof of \cite[Proposition~5.3]{topological-branching}. Let $Z\subset |L|_{geom}$ and $\mathcal{P}^-_\alpha(L,\epsilon,Z)_{co} = \mathcal{P}^-_\alpha(L,\epsilon) \cap \ttop_{co}([0,1],Z)$ equipped with the relative topology for the inclusion $\mathcal{P}^-_\alpha(L,\epsilon,Z)_{co}\subset \ttop_{co}([0,1],Z)$. Let $\mathcal{P}^-_\alpha(L,\epsilon,Z)=k_\Delta(\mathcal{P}^-_\alpha(L,\epsilon,Z)_{co})$. The point is that $\mathcal{P}^-_\alpha(L,\epsilon,|K|_{geom})=\mathcal{P}^-_\alpha(K,\epsilon)$. From the pushout diagram of spaces 
\[\begin{tikzcd}[row sep=3em,column sep=3em]
	{|\de\square[n]|_{geom}} \arrow[r,"|g|_{geom}"] \arrow[d,"\subset"'] & {|K|_{geom}} \arrow[d] \\
	{|\square[n]|_{geom}} \arrow[r,"|\widehat{g}|_{geom}"] & \cocartesian |L|_{geom}
\end{tikzcd}\]
we deduce that the map $|K|_{geom}\to |L|_{geom}$ is a q-cofibration of spaces (see e.g. \cite[Proposition~3.12]{DirectedDegeneracy} for further explanations). Therefore it is a h-cofibration. By \cite[Theorem~2]{vstrom1}, there exists an open subset $U$ of $|L|_{geom}$ such that $|K|_{geom}\subset U$ and such that $|K|_{geom}$ is a retract of $U$. We then follow the proof of \cite[Proposition~5.3]{topological-branching}. 
\begin{enumerate}[leftmargin=*]
	\item The space $\mathcal{P}^-_\alpha(K,\epsilon)$ being a retract of the space $\mathcal{P}^-_\alpha(L,\epsilon,U)$, there is a closed inclusion $\mathcal{P}^-_\alpha(K,\epsilon) \to \mathcal{P}^-_\alpha(L,\epsilon,U)$. For the same reason, the map $|K|_{geom}\to U$ is a closed inclusion as well.
	\item Let $(\gamma_n)_{n\geq 0}$ be a sequence of $\mathcal{P}^-_\alpha(K,\epsilon)$ which converges to $\gamma_\infty$ in $\mathcal{P}^-_\alpha(L,\epsilon)$. Since $|K|_{geom}$ is closed in $|L|_{geom}$ which is Hausdorff, and since $(\gamma_n)_{n\geq 0}$ converges pointwise to $\gamma_\infty\in \ttop([0,1],|L|_{geom})$, one has $\gamma_\infty([0,1])\subset |K|_{geom}$. Since $|K|_{geom}\subset U$ is a closed inclusion and since $U$ is open in $|L|_{geom}$, the composite map $|K|_{geom}\to U \to |L|_{geom}$ is a $\Delta$-inclusion, $U$ being equipped with the relative topology. This implies that $\gamma_\infty\in\ttop([0,1],|K|_{geom})$. Since $\mathcal{P}^-_\alpha(L,\epsilon,U)_{co}$ is an open subset of $\mathcal{P}^-_\alpha(L,\epsilon)_{co}$ by definition of the compact-open topology, $\mathcal{P}^-_\alpha(L,\epsilon,U)$ is an open subspace of $\mathcal{P}^-_\alpha(L,\epsilon)$ by Lemma~\ref{lem:open-delta-k}. Consequently, there exists $N\geq 0$ such that $\gamma_n\in \mathcal{P}^-_\alpha(L,\epsilon,U)$ for all $n\geq N$ and, moreover, the sequence $(\gamma_n)_{n\geq N}$ is convergent in $\mathcal{P}^-_\alpha(L,\epsilon,U)$.
	\item Thus the sequence $(\gamma_n)_{n\geq 0}$ converges to $\gamma_\infty$ in $\mathcal{P}^-_\alpha(K,\epsilon)$. Since $|L|_{geom}$ is Hausdorff, it is $\Delta$-Hausdorff. We deduce that $\mathcal{P}^-_\alpha(L,\epsilon)$ is $\Delta$-Hausdorff as well. The proof is then complete thanks to Theorem~\ref{thm:recognizing-closed-inclusions-1}.
\end{enumerate}
\epf

\begin{prop} \label{prop:compacity}
	The set $N_0(\epsilon)$ is empty. Let $n\geq 1$. The set $N_n(\epsilon)=\mathcal{P}^-_{0_n}(\square[n],\epsilon)$ equipped with the compact-open topology is $\Delta$-generated, $\Delta$-Hausdorff, metrizable, contractible, compact and sequentially compact. 
\end{prop}

\bpf
The compact-open topology is metrizable with the distance of the uniform convergence by \cite[Proposition~A.13]{MR1867354}. Therefore it is first countable. Consider a ball $B(\gamma,\eta)$ for this metric. Let $\gamma'\in B(\gamma,\eta)$. Then each convex combination $(1-u)\gamma+u\gamma'$ of directed paths of $N_n(\epsilon)$ belongs to $N_n(\epsilon)$. Moreover, for all $t\in [0,\epsilon]$ and all $i\in \{1,\dots,n\}$, one has 
\[
|((1-u)\gamma_i+u\gamma'_i)(t) - \gamma_i(t)| = u|\gamma'_i(t) - \gamma_i(t)| < u\eta \leq \eta.
\]
This means that the space $N_n(\epsilon)$ is locally path-connected. By \cite[Proposition~3.11]{MR3270173}, it is $\Delta$-generated, and also $\Delta$-Hausdorff, being metrizable. Each directed path of $N_n(\epsilon)$ can be extended using a straight line to a natural directed path from $[0,n]$ to $[0,1]^n$ ending to $1_n$. We obtain a continuous map $N_n(\epsilon)\to N_n$, where $N_n$ is the space of natural directed paths of $[0,1]^n$ going from $0_n$ to $1_n$. By restricting to $[0,\epsilon]$, we obtain a continuous map $N_n\to N_n(\epsilon)$ such that the composite map $N_n(\epsilon)\to N_n\to N_n(\epsilon)$ is the identity. Thus $N_n(\epsilon)$ is a retract of $N_n$. This implies that $N_n(\epsilon)$ is contractible, $N_n$ being contractible by \cite[Proposition~4.14]{NaturalRealization}. This also implies that $N_n(\epsilon)$ is sequentially compact, $N_n$ being sequentially compact by \cite[Proposition~4.14]{NaturalRealization}. Being metrizable, we deduce that $N_n(\epsilon)$ is compact.  
\epf

\begin{prop} \label{prop:closedness}
	The subset $\mathcal{P}^-_{0_n}(\de\square[n],\epsilon)$ equipped with the relative topology is a closed subset of $\mathcal{P}^-_{0_n}(\square[n],\epsilon)$ which is $\Delta$-generated, $\Delta$-Hausdorff, metrizable, compact and sequentially compact.
\end{prop}

\bpf
By Proposition~\ref{prop:P-pushout}, the inclusion $\mathcal{P}^-_{0_n}(\de\square[n],\epsilon) \subset \mathcal{P}^-_{0_n}(\square[n],\epsilon)$ is a closed inclusion of $\Delta$-generated spaces. Thus, the relative topology on $\mathcal{P}^-_{0_n}(\de\square[n],\epsilon)$ with respect to this inclusion is $\Delta$-generated. The rest of the statement is clear.
\epf

\section{Categorical properties of the homotopy branching space}
\label{section:cat-branching}

We recall the useful lemma: 

\begin{lem} \label{lem:converge} \cite[Lemma~2]{Moore3}
	Let $X$ be a sequential topological space. Let $x_\infty\in X$. Let $(x_n)_{n\geq 0}$ be a sequence such that $x_\infty$ is a limit point of all subsequences. Then the sequence $(x_n)_{n\geq 0}$ converges to $x_\infty$.
\end{lem}

\begin{prop} \label{prop:G-pushout-2}
	Let $n\geq 0$. Consider a pushout diagram of precubical sets
	\[
	\begin{tikzcd}[row sep=3em,column sep=3em]
		\de\square[n] \arrow[r,"g"] \arrow[d] & K \arrow[d,"f"] \\
		\square[n] \arrow[r,"\widehat{g}"] & \cocartesian L
	\end{tikzcd}
	\]
	Then there is a pushout diagram of spaces 
	\[
	\begin{tikzcd}[row sep=3em,column sep=3em]
		\mathcal{P}^-_{0_n}(\de\square[n],\epsilon) \arrow[r] \arrow[d] & \mathcal{P}^-_{g(0_n)}(K,\epsilon) \arrow[d] \\
		\mathcal{P}^-_{0_n}(\square[n],\epsilon) \arrow[r] & \cocartesian \mathcal{P}^-_{\widehat{g}(0_n)}(L,\epsilon)
	\end{tikzcd}
	\]
\end{prop}

\bpf There is nothing to prove for $n=0$. Assume that $n\geq 1$. There is the commutative diagram of topological spaces 
\[
\begin{tikzcd}[row sep=3em,column sep=3em]
	\mathcal{P}^-_{0_n}(\de\square[n],\epsilon) \arrow[r] \arrow[d] & \mathcal{P}^-_{g(0_n)}(K,\epsilon) \arrow[d,"h_1"] \arrow[ddr,bend left=20pt,"k_1"]\\
	\mathcal{P}^-_{0_n}(\square[n],\epsilon) \arrow[r,"h_2"] \arrow[rrd,bend right=20pt,"k_2"']&  \cocartesian Z \arrow[rd,dashed,"\exists ! k"]\\ 
	&& \mathcal{P}^-_{\widehat{g}(0_n)}(L,\epsilon)
\end{tikzcd}
\]
Let $[c]_{geom}\phi \in \mathcal{P}^-_{\widehat{g}(0_n)}(L,\epsilon)$. Then $c\in \mathcal{C}^-_{\widehat{g}(0_n)}(L) = \mathcal{C}^-_{g(0_n)}(K) \cup \{\widehat{g}\}$. Thus the map $k$ is surjective. Let $z_1,z_2\in Z$ such that $k(z_1)=k(z_2)$. There are three possible cases:
\begin{enumerate}
	\item $h_1(\overline{z}_i)=z_i$ for $i=1,2$; then $k_1(\overline{z}_1)=k_1(\overline{z}_2)$; the map $k_1$ being one-to-one by Proposition~\ref{prop:P-pushout}, we deduce that $\overline{z}_1=\overline{z}_2$ and consequently that $z_1=z_2$.
	\item $h_2(\phi_i)=z_i$ for $i=1,2$; in this case, one has $k(z_i)=|\widehat{g}|_{geom}\phi_i$ for $i=1,2$; we deduce that $\phi_1=\phi_2$ and therefore that $z_1=z_2$.
	\item $h_1(|c|_{geom}\phi_1)=z_1$ and $h_2(\phi_2)=z_2$; we deduce that $|c|_{geom}\phi_1=|\widehat{g}|_{geom}\phi_2$; this implies that the image of $\phi_2:[0,\epsilon]\to [0,1]^n$ is included in the boundary of $[0,1]^n$; we then deduce that $z_1=z_2$ like in the first case.
\end{enumerate} 
We have proved that the map $k:Z\to \mathcal{P}^-_{\widehat{g}(0_n)}(L,\epsilon)$ is a continuous bijection. Consider a sequence $(|c_k|_{geom}\phi_k)_{k\geq 0}$ of $\mathcal{P}^-_{\widehat{g}(0_n)}(L,\epsilon)$ converging to $\gamma_\infty$. There are two mutually exclusive cases: 
\begin{enumerate}
	\item $\{k\geq 0\mid c_k=\widehat{g}\}$ infinite: by extracting a subsequence, we can suppose that the sequence is of the form $(|\widehat{g}|_{geom}\phi_k)_{k\geq 0}$ where $\phi_k\in \mathcal{P}^-_{0_n}(\square[n],\epsilon)$. By Proposition~\ref{prop:compacity}, we can suppose that the sequence $(\phi_k)_{k\geq 0}$ has a limit. We deduce that $(|\widehat{g}|_{geom}\phi_k)_{k\geq 0}$ is convergent as well. Since $k$ is bijective, this limit is necessarily $\gamma_\infty$.
	\item $\{k\geq 0\mid c_k=\widehat{g}\}$ finite: by extracting a subsequence, we can suppose that the sequence $(|c_k|_{geom}\phi_k)_{k\geq 0}$ belongs to $k_1(\mathcal{P}^-_{g(0_n)}(K,\epsilon))$; the map $k_1$ being a closed inclusion by Proposition~\ref{prop:P-pushout}, the sequence $(|c_k|_{geom}\phi_k)_{k\geq 0}$ converges to $\gamma_\infty$ in $\mathcal{P}^-_{g(0)}(K,\epsilon)$; thus it converges in $Z$ as well; we have prove that the sequence $(|c_k|_{geom}\phi_k)_{k\geq 0}$ of $Z$ has a limit point which is necessarily $\gamma_\infty$. 
\end{enumerate}
We have proved that every subsequence of $(|c_k|_{geom}\phi_k)_{k\geq 0}$ of $Z$ as a limit point which is necessarily $\gamma_\infty$. By Lemma~\ref{lem:converge}, $Z$ being sequential, the sequence $(|c_k|_{geom}\phi_k)_{k\geq 0}$ converges to $\gamma_\infty$ in $Z$. We deduce that $Z$ and $\mathcal{P}^-_{\widehat{g}(0)}(L,\epsilon)$  have the same convergent sequences. Since both spaces are sequential, being $\Delta$-generated, we conclude that the continuous bijection $k:Z\to \mathcal{P}^-_{\widehat{g}(0)}(L,\epsilon)$ is a homeomorphism.
\epf

\begin{cor} \label{cor:P-pushout}
	Let $n\geq 0$. Consider a pushout diagram of precubical sets
	\[
	\begin{tikzcd}[row sep=3em,column sep=3em]
		\de\square[n] \arrow[r,"g"] \arrow[d] & K \arrow[d] \\
		\square[n] \arrow[r,"\widehat{g}"] & \cocartesian L
	\end{tikzcd}
	\]
	Then there is a pushout diagram of spaces 
	\[
	\begin{tikzcd}[row sep=3em,column sep=3em]
		\mathcal{P}^-(\de\square[n],\epsilon) \arrow[r] \arrow[d] & \mathcal{P}^-(K,\epsilon) \arrow[d] \\
		\mathcal{P}^-(\square[n],\epsilon) \arrow[r] & \cocartesian \mathcal{P}^-(L,\epsilon)
	\end{tikzcd}
	\]
\end{cor}

\bpf
By Proposition~\ref{prop:P-pushout}, there are the equalities $\mathcal{P}^-_\alpha(\de\square[n],\epsilon)=\mathcal{P}^-_\alpha(\square[n],\epsilon)$ and $\mathcal{P}^-_{g(\alpha)}(K,\epsilon) = \mathcal{P}^-_{\widehat{g}(\alpha)}(L,\epsilon)$ for $a\neq 0_n$. The proof is complete thanks to Proposition~\ref{prop:G-pushout-2}.
\epf

\begin{prop} \label{prop:P-colim}
	Let $\widetilde{X}:\lambda\to\square^{op}\set$ be a transfinite tower (i.e. a colimit-preserving tower) of precubical sets such that each map $\widetilde{X}_\nu\to \widetilde{X}_{\nu+1}$ is a pushout of the map of precubical sets $\de\square[n_\nu]\subset \square[n_\nu]$ for some $n_\nu\geq 0$. Let $\nu$ be a limit ordinal with $\nu\leq \lambda$. Then the canonical continuous map 
	\[
	\Phi_\nu:\liminj_{\mu<\nu} \mathcal{P}^-(\widetilde{X}_\mu,\epsilon) \longrightarrow \mathcal{P}^-(\widetilde{X}_\nu,\epsilon)
	\]
	is a homeomorphism.
\end{prop}

\bpf
There is the equality of precubical sets $\widetilde{X}_\nu=\bigcup_{\mu<\nu} \widetilde{X}_\mu$. This implies that $\Phi_\nu$ is surjective. Let $|c_i|_{geom}\phi_i$ in $\mathcal{P}^-(\widetilde{X}_{\mu_i},\epsilon)$ with $i=1,2$ for $\mu_1,\mu_2<\nu$ such that $\Phi_\nu(|c_1|_{geom}\phi_1)=\phi_\nu(|c_2|_{geom}\phi_2)$. By taking $\max(\mu_1,\mu_2)<\nu$, we can suppose that $\mu_1=\mu_2=\mu$. Then the two continuous maps $|c_1|_{geom}\phi_1,|c_2|_{geom}\phi_2:[0,1]\to |\widetilde{X}_\mu|_{geom}$ induces the same continuous map $[0,1]\to |\widetilde{X}_\mu|_{geom} \to |\widetilde{X}_\nu|_{geom}$. Since the map $|\widetilde{X}_\mu|_{geom} \to |\widetilde{X}_\nu|_{geom}$ is a q-cofibration (see e.g. \cite[Proposition~3.12]{DirectedDegeneracy} for further explanations), it is one-to-one. This implies that $|c_1|_{geom}\phi_1)=|c_2|_{geom}\phi_2$ as continuous maps from $[0,1]$ to $|\widetilde{X}_\mu|_{geom}$. Thus $\Phi_\nu$ is one-to-one. We have proved that $\Phi_\nu$ is a continuous bijection. Consider a set map $f:[0,1]\to \liminj_{\mu<\nu} \mathcal{P}^-(\widetilde{X}_\mu,\epsilon)$ such that the composite set map $\Phi_\nu f$ is continuous. It gives rise to a continuous map $\widetilde{f}:[0,1]\p [0,1]\to |\widetilde{X}_\nu|_{geom}$. Since the tower $|\widetilde{X}|_{geom}$ consists of q-cofibrations which are closed $T_1$-inclusions, the map $\widetilde{f}:[0,1]\p [0,1]\to |\widetilde{X}_\nu|_{geom}$ factors as a composite $\widetilde{f}:[0,1]\p [0,1]\to |\widetilde{X}_\mu|_{geom}\to|\widetilde{X}_\nu|_{geom}$ for some $\mu<\nu$ by \cite[Proposition~2.4.2]{MR99h:55031}, $[0,1]\p [0,1]$ being compact. By adjunction, we obtain that $\Phi_\nu f$ factors as a composite of continuous maps $[0,1]\to \mathcal{P}^-(\widetilde{X}_\nu,\epsilon)\cap \ttop([0,1],|\widetilde{X}_\mu|_{geom})= \mathcal{P}^-(\widetilde{X}_\mu,\epsilon) \to \mathcal{P}^-(\widetilde{X}_\nu,\epsilon)$. Hence $f$ is continuous. Therefore, the continuous bijection $\Phi_\nu$ is a $\Delta$-inclusion. Thus, it is a homeomorphism.
\epf

\bth \label{thm:G-colimit-preserving}
The functor $\mathcal{P}^-(-,\epsilon):\square^{op}\set \to \top$ is colimit-preserving.
\eth

\bpf
The cocubical object $\mathcal{P}^-(\square[*],\epsilon)$ gives rise to a colimit-preserving functor 
\[
\widehat{\mathcal{P}}(K,\epsilon)=\int^{[n]\in \square} K_n.\mathcal{P}^-(\square[n],\epsilon).
\]
Using Proposition~\ref{prop:P-colim} and Corollary~\ref{cor:P-pushout}, we prove by an easy induction that there is a homeomorphism $\widehat{\mathcal{P}}(K,\epsilon)\iso \mathcal{P}^-(K,\epsilon)$ for all precubical sets $K$. Hence the proof is complete.
\epf

\section{Homotopical properties of the homotopy branching space}
\label{section:homotopic-branching}

\begin{prop}\label{prop:cubical-cofibration-1}
	Let $n\geq 0$. The map $\mathcal{P}^-_{0_n}(\de\square[n],\epsilon) \to \mathcal{P}^-_{0_n}(\square[n],\epsilon)$ is an m-cofibration between m-cofibrant spaces.
\end{prop}

\bpf For $n=0$, there is nothing to prove. Assume that $n\geq 1$. Consider the locally convex topological vector space $E=\top([0,\epsilon],\mathbb{R}^n)$ of continuous maps from $[0,\epsilon]$ to $\mathbb{R}^n$  equipped with the sup norm (i.e. the underlying topology is the compact-open topology). Consider $[0,1]^n$ as a subset of $\mathbb{R}^n$. As already noticed in the proof of Proposition~\ref{prop:compacity}, the nonempty space $N_n(\epsilon)=\mathcal{P}^-_{0_n}(\square[n],\epsilon)$ is closed by convex sum. Thus it is an AR by the Dugundji extension theorem \cite[Theorem~6.1.1]{Sakai}. Since $0<\epsilon<1$, one has 
\[
\mathcal{P}^-_{0_n}(\de\square[n],\epsilon) = \bigcup_{i=1}^n F_i
\]
where $F_i\subset \mathcal{P}^-_{0_n}(\de\square[n],\epsilon)$ consists of the natural directed paths $\gamma=(\gamma_1,\dots,\gamma_n)$ such that $\gamma_i=0$. Each finite intersection of sets of $\{F_1,\dots,F_n\}$ is empty, and then an ANR, or nonempty convex in $E$, and then an AR by \cite[Theorem~6.1.1]{Sakai} again. By \cite[6.2.10 (5)]{Sakai}, any AR being an ANR, we deduce that $\mathcal{P}^-_{0_n}(\de\square[n],\epsilon)$ is an ANR. Using Proposition~\ref{prop:closedness}, we have proved that $\mathcal{P}^-_{0_n}(\de\square[n],\epsilon)$ is a closed ANR and that $\mathcal{P}^-_{0_n}(\square[n],\epsilon)$ is an ANR of the metrizable space $E$. Thus the pair $(\mathcal{P}^-_{0_n}(\square[n],\epsilon),\mathcal{P}^-_{0_n}(\de\square[n],\epsilon))$ is an NDR-pair by \cite[Theorem~4.2.15]{AlgTopo-Homotopical-PoV}. By \cite[Theorem~2]{vstrom1}, we deduce that the map $\mathcal{P}^-_{0_n}(\de\square[n],\epsilon) \to \mathcal{P}^-_{0_n}(\square[n],\epsilon)$ is a h-cofibration. By \cite[Corollary~6.6.5]{Sakai}, the topological spaces $\mathcal{P}^-_{0_n}(\de\square[n],\epsilon)$ and  $\mathcal{P}^-_{0_n}(\square[n],\epsilon)$ are m-cofibrant, being ANR. Thanks to \cite[Corollary~3.12]{mixed-cole}, we deduce that the h-cofibration $\mathcal{P}^-_{0_n}(\de\square[n],\epsilon) \to \mathcal{P}^-_{0_n}(\square[n],\epsilon)$ is an m-cofibration. 
\epf

\begin{prop} \label{prop:cubical-cofibration-2}
	Let $n\geq 0$. Consider a pushout diagram of precubical sets
	\[
	\begin{tikzcd}[row sep=3em,column sep=3em]
		\de\square[n] \arrow[r,"g"] \arrow[d] & K \arrow[d] \\
		\square[n] \arrow[r,"\widehat{g}"] & \cocartesian L
	\end{tikzcd}
	\]
	Then the vertical maps of the pushout diagram of spaces 
	\[
	\begin{tikzcd}[row sep=3em,column sep=3em]
		\mathcal{P}^-(\de\square[n],\epsilon) \arrow[r] \arrow[d,rightarrowtail] & \mathcal{P}^-(K,\epsilon) \arrow[d,rightarrowtail] \\
		\mathcal{P}^-(\square[n],\epsilon) \arrow[r] & \cocartesian \mathcal{P}^-(L,\epsilon)
	\end{tikzcd}
	\]
	are m-cofibrations. Consequently, for all precubical sets $K$, the space $\mathcal{P}^-(K,\epsilon)$ is m-cofibrant.
\end{prop}

\bpf
For $\alpha\in \{0,1\}^n\backslash \{0_n\}$ and $n\geq 1$, there is a homeomorphism $\mathcal{P}^-_{\alpha}(\de\square[n],\epsilon) = \mathcal{P}^-_{\alpha}(\square[n],\epsilon)$ by Proposition~\ref{prop:P-pushout}. Using Proposition~\ref{prop:cubical-cofibration-1}, we deduce that the left vertical map is an m-cofibration. Hence the right vertical map is an m-cofibration as well. The last assertion is obtained by an easy induction.
\epf

\begin{cor} \label{cor:cubical-cofibration-3}
	Let $K\to L$ be a transfinite composition of pushouts of maps of $\{\de\square[n]\subset \square[n]\mid n\geq 0\}$. Then the map $\mathcal{P}^-(K,\epsilon)\to \mathcal{P}^-(L,\epsilon)$ is an m-cofibration of spaces.
\end{cor}

\begin{nota} \label{nota:degenerate-branching}
	The cocubical space $\mathbb{B}:\square\to \top$ is defined as follows:
	\begin{itemize}
		\item $\mathbb{B}([0]) = \varnothing$;
		\item for $n\geq 1$, $\mathbb{B}([n]) =\{0,1\}^n\backslash \{1_n\}$
		\item for $f:[m]\to [n]$, $\mathbb{B}(f)(\alpha) = f(\alpha)$.
	\end{itemize}
\end{nota}

\bth \label{thm:P-homotopy-invariant}
Let $\epsilon,\epsilon'\in ]0,1[$. For all precubical sets $K$, there is a natural homotopy equivalence between m-cofibrant spaces 
\[
\mathcal{P}^-(K,\epsilon) \simeq \mathcal{P}^-(K,\epsilon').
\]
\eth

\bpf
By Theorem~\ref{thm:tout} and Proposition~\ref{prop:cubical-cofibration-2}, the two cocubical spaces $\mathcal{P}^-(\square[*],\epsilon)$ and $\mathcal{P}^-(\square[*],\epsilon')$ are Reedy cofibrant. By Proposition~\ref{prop:compacity}, the space $\mathcal{P}^-_{0_n}(\square[n],\epsilon)$ is contractible for all $n\geq 1$. More generally, for $n\geq 1$ and $\alpha\in \mathbb{B}([n])$, there is the homeomorphism $\mathcal{P}^-_{\alpha}(\square[n],\epsilon) \iso \mathcal{P}^-_{0_d}(\square[d],\epsilon)$ with $d=(1-\alpha_1)+\dots +(1-\alpha_n)$ with $\alpha=(\alpha_1,\dots,\alpha_n)$. Thus by Proposition~\ref{prop:compacity}, there is a map of cocubical spaces $\mathcal{P}^-(\square[*],\epsilon) \to \mathbb{B}$ which is an objectwise homotopy equivalence. Since all spaces $\mathbb{B}([n])$ are discrete, we obtain the following diagram of Reedy trivial fibrations of cocubical spaces
\[
\begin{tikzcd}[row sep=small,column sep=3em]
	\mathcal{P}^-(\square[*],\epsilon) \arrow[r,twoheadrightarrow,"\simeq"] & \mathbb{B}([*]) & \mathcal{P}^-(\square[*],\epsilon') \arrow[l,twoheadrightarrow,"\simeq"']
\end{tikzcd}
\]
The proof is complete thanks to Theorem~\ref{thm:tout}.
\epf

There is a tower of topological spaces 
\[
\begin{tikzcd}[row sep=small,column sep=3em]
	\mathcal{P}^-(K,1/2)\arrow[r] & \mathcal{P}^-(K,1/3) \arrow[r]& \mathcal{P}^-(K,1/4) \arrow[r]& \dots
\end{tikzcd}
\]
induced by the restriction of natural directed paths.

\begin{prop} \label{prop:ho-germ}
	For all precubical sets $K$, there is the homeomorphism 
	\[
	\mathcal{G}^-(K) \iso \liminj_{n\geq 0} \mathcal{P}^-(K,1/(n+2))
	\]
	where $\mathcal{G}^-(K)$ is the quotient of $\mathcal{P}^-(K,1/2)$ by the equivalence relation $\simeq^-$ defined in Proposition~\ref{prop:no-germs}. Moreover there is the homotopy equivalence 
	\[
	\mathcal{P}^-(K,1/2) \simeq \holiminj_{n\geq 0} \mathcal{P}^-(K,1/(n+2))
	\]
\end{prop}

\bpf
The homeomorphism is a straightforward consequence of the definition of the equivalence relation $\simeq^-$. By Theorem~\ref{thm:P-homotopy-invariant}, the homotopy colimit in the statement of the proposition is the homotopy colimit of the constant diagram $D:(\mathbb{N},\leq)\to \top$ defined by $D(n)=\mathcal{P}^-(K,1/2)$. By \cite[Proposition~18.1.6]{ref_model2}, we then have ($B$ meaning the classifying space) \[\holiminj_{n\geq 0} \mathcal{P}^-(K,1/(n+2))\simeq B(\mathbb{N},\leq)^{op} \p \mathcal{P}^-(K,1/2).\] The small category $(\mathbb{N},\leq)$ is contractible, having an initial object. Hence the proof is complete.
\epf

The homotopy branching space of a precubical set $K$ can therefore be interpreted as the space of homotopy germs of short natural directed paths.

\section{Homotopy branching space of a flow}
\label{section:flow-branching}

\bd \cite[Definition~4.11]{model3} \label{def:flow}
A \textit{flow} $X$ is a small topologically enriched semicategory. Its set of objects (preferably called \textit{states}) is denoted by $X^0$ and the space of morphisms (preferably called \textit{execution paths}) from $\alpha$ to $\beta$ is denoted by $\P_{\alpha,\beta}X$. For any $x\in \P_{\alpha,\beta}X$, let $s(x)=\alpha$ and $t(x)=\beta$. The category is denoted by $\dtop$. 
\ed

The category $\dtop$ is locally presentable. Every set can be viewed as a flow with an empty path space. This give rise to a functor from sets to flows which is limit-preserving and colimit-preserving. More generally, any poset can be viewed as a flow, with a unique execution path from $u$ to $v$ if and only if $u<v$. This gives rise to a functor from the category of posets together with the strictly increasing maps to flows.

\begin{nota} \label{nota:glob}
	For any topological space $Z$, the flow $\glob(Z)$ is the flow having two states $0$ and $1$ and such that the only nonempty space of execution paths, when $Z$ is nonempty, is $\P_{0,1}\glob(Z)=Z$. It is called \textit{the globe of $Z$}. Let $\vI = \glob(\mathbf{D}^0)$.
\end{nota}

We need to recall:

\bth \label{thm:three} \cite[Theorem~2.4]{NaturalRealization} Let $r\in \{q,m,h\}$. Then there exists a unique model structure on $\dtop$ such that: 
\begin{itemize}
	\item A map of flows $f:X\to Y$ is a weak equivalence if and only if $f^0:X^0\to Y^0$ is a bijection and for all $(\alpha,\beta)\in X^0\p X^0$, the continuous map $\P_{\alpha,\beta}X\to \P_{f(\alpha),f(\beta)}Y$ is a weak equivalence of the r-model structure of $\top$.
	\item A map of flows $f:X\to Y$ is a fibration if and only if for all $(\alpha,\beta)\in X^0\p X^0$, the continuous map $\P_{\alpha,\beta}X\to \P_{f(\alpha),f(\beta)}Y$ is a fibration of the r-model structure of $\top$.
\end{itemize}
This model structure is accessible and all objects are fibrant. Moreover, this model structure is simplicial. It is called the r-model structure of $\dtop$.
\eth

\bth \label{thm:def-Cminus} 
Let $r\in \{q,m,h\}$.  Let $Z$ be a topological space. The data $C^-(Z)^0=\{0\}$, $\P_{0,0}Z=Z$ and $x*y=x$ assemble to a flow $C^-(Z)$. The mapping $Z\mapsto C^-(Z)$ gives rise to a right Quillen adjoint $C^-:\top\to \dtop$ from the r-model structure of $\top$ to the r-model structure of $\dtop$. 
\eth

\bpf
It suffices to prove that the functor $C^-$ takes (trivial) r-fibrations of spaces to (trivial) r-fibrations of flows. For r=q, the details are explained in \cite[Theorem~5.5]{exbranching}.
\epf

The left Quillen adjoint $\P^-:\dtop \to \top$ is defined by taking a flow $X$ to the quotient space 
\[
\P^- X = \coprod_{\alpha\in X^0} \P^-_\alpha X
\]
where $\P_\alpha X$ is the quotient of the space $\coprod_{\beta\in X^0} \P_{\alpha,\beta}X$ by the equivalence relation generated by the identifications $u*v=u$. By \cite[Theorem~4.1]{exbranching}, the branching space functor $\P^-:\dtop \to \top$ is badly behaved with respect to the weak equivalences. Indeed, there exists a weak equivalence of flows $f:X\to Y$ such that $\P^-f:\P^-X\to \P^-Y$ is not a weak homotopy equivalence of spaces. Hence the following definition.

\bd \label{def:branching-flow}
Let $X$ be a flow. The topological space $\P^-X$ is called the \textit{branching space} of $X$. The topological space $\hop^-X=\P^- X^{cof}$ is called the \textit{homotopy branching space} of $X$ where $(-)^{cof}$ is some q-cofibrant replacement of $X$. The latter space is unique only up to homotopy equivalence.
\ed

\bd \label{def-rea-flow} \cite[Definition~3.6]{NaturalRealization} Let $r\in \{q,m,h\}$. A functor $F:\square^{op}\set \to \dtop$ is an \textit{r-realization functor} if it satisfies the following properties:
\begin{itemize}
	\item $F$ is colimit-preserving.
	\item For all $n\geq 0$, the map of flows $F(\de\square[n])\to F(\square[n])$ is an r-cofibration.
	\item There is an objectwise weak equivalence of cocubical flows $F(\square[*])\to \{0<1\}^*$ in the r-model structure of $\dtop$.
\end{itemize}
Every flow $F(K)$ is r-cofibrant.
\ed

\bth \label{thm:top-flow-homotopy-equivalent}
Let $F:\square^{op}\set\to \dtop$ be a q-realization functor. Then for all precubical sets $K$, there is a natural homotopy equivalence 
\[
\mathcal{P}^-(K,\epsilon) \simeq \P^-F(K).
\]
\eth

\bpf
By definition of a q-realization functor, for all $n\geq 0$, the map of flows $F(\de\square[n])\to F(\square[n])$ is a q-cofibration. Thus, $\P^-$ being a left Quillen adjoint by Theorem~\ref{thm:def-Cminus}, the map $\P^-F(\de\square[n])\to \P^-F(\square[n])$ is a q-cofibration of spaces for all $n\geq 0$. We deduce that the cocubical space $\P^-F(\square[*])$ is Reedy q-cofibrant by Theorem~\ref{thm:tout}~(1). The point is that the poset $\{0<1\}^n$ is finite bounded. By \cite[Theorem~9.3]{3eme}, the space $\P^-_\alpha F(\{0<1\}^n)$ is contractible for $\alpha\in \{0,1\}^n \backslash\{1_n\}$ and empty for $\alpha=1_n$ when $n\geq 1$. Thus there is a trivial Reedy fibration of cocubical spaces $\P^- F(\square[*])\to \mathbb{B}([*])$. We conclude using Theorem~\ref{thm:tout}~(3).
\epf

\begin{cor} \label{cor:top-flow-homotopy-equivalent}
	Let $r\in \{q,m,h\}$. Let $F:\square^{op}\set\to \dtop$ be an r-realization functor. Then for all precubical sets $K$, there is a natural homotopy equivalence 
	\[
	\mathcal{P}^-(K,\epsilon) \simeq \P^-F(K).
	\]
\end{cor}

\bpf
Choose a q-realization functor $|-|_q$. Then it is an r-realization functor, any q-cofibration being an r-cofibration. By \cite[Theorem~3.8]{NaturalRealization}, there is a natural homotopy equivalence of flows $|K|_q\simeq F(K)$. Since $|K|_q$ is q-cofibrant, it is r-cofibrant, as well as $F(K)$. Thus there is a natural weak equivalence $\P^-|K|_q \to \P^-F(K)$, $\P^-$ being a left Quillen adjoint. Being also r-fibrant, the natural weak equivalence is a natural homotopy equivalence. The proof is complete thanks to Theorem~\ref{thm:top-flow-homotopy-equivalent}.
\epf

\begin{cor} \label{cor:hop-not-necessary}
	Let $r\in \{q,m,h\}$. Let $F:\square^{op}\set\to \dtop$ be an r-realization functor. There is a natural homotopy equivalence $\hop^-F(K)\simeq \P^-F(K)$ for any precubical set $K$ and any cofibrant replacement functor $(-)^{cof}$ of the q-model category of flows.
\end{cor}

\bpf
By Corollary~\ref{cor:top-flow-homotopy-equivalent}, it suffices to prove the statement for a q-realization functor $F:\square^{op}\set\to \dtop$. Then the natural map $F(K)^{cof}\to F(K)$ is a weak equivalence between two q-cofibrant flows. The functor $\P^-:\dtop\to \top$ being a left Quillen functor, the map $\hop^-F(K)\to \P^-F(K)$ is a weak homotopy equivalence between q-cofibrant spaces. Hence it is a homotopy equivalence.
\epf

\section{Homotopy branching space and cubical subdivision}
\label{section:sub}

\begin{figure}
	\def\n{5}
	\begin{tikzpicture}[black,scale=4,pn/.style={circle,inner sep=0pt,minimum width=4pt,fill=dark-red}]
		\fill [color=gray!15] (0,0) -- (0,1) -- (1,1) -- (1,0) -- cycle;
		\draw[dark-red][-] (0,0) -- (0,1/3) -- (1/3,1/3) -- (1/3,0) -- cycle;
		\draw[dark-red][-] (1/3,0) -- (1/3,1/3) -- (2/3,1/3) -- (2/3,0) -- cycle;
		\draw[dark-red][-] (2/3,0) -- (2/3,1/3) -- (1,1/3) -- (1,0) -- cycle;
		\draw[dark-red][-] (0,1/3) -- (0,2/3) -- (1/3,2/3) -- (1/3,1/3) -- cycle;
		\draw[dark-red][-] (1/3,1/3) -- (1/3,2/3) -- (2/3,2/3) -- (2/3,1/3) -- cycle;
		\draw[dark-red][-] (2/3,1/3) -- (2/3,2/3) -- (1,2/3) -- (1,1/3) -- cycle;
		\draw[dark-red][-] (0,2/3) -- (0,1) -- (1/3,1) -- (1/3,2/3) -- cycle;
		\draw[dark-red][-] (1/3,2/3) -- (1/3,1) -- (2/3,1) -- (2/3,2/3) -- cycle;
		\draw[dark-red][-] (2/3,2/3) -- (2/3,1) -- (1,1) -- (1,2/3) -- cycle;
	\end{tikzpicture}
	\caption{Cubical subdivision $\Sub_3(\square[2])$}
	\label{fig:sub3-2}
\end{figure}

Let $n\geq 1$ and $p\geq 1$. Consider the equality of spaces
\[
[0,p]^n = \bigcup_{\substack{a_{i_1},\dots,a_{i_n}\in \{0,\dots,p-1\}}} [a_{i_1},a_{i_1}+1]\p \dots \p [a_{i_n},a_{i_n}+1]
\]
Then consider the precubical set $\Sub_p(\square[n])$ such that the set of $p^n$ $n$-cubes is $[a_{i_1},a_{i_1}+1]\p \dots \p [a_{i_n},a_{i_n}+1]$ for $a_{i_1},\dots,a_{i_n}\in \{0,\dots,p-1\}$ and such that the face maps are obtained by replacing segments $[a,b]$ of the products by $a$ or $b$. In particular, one has $\Sub_p(\square[n])_0=\{0,1,\dots,p\}^n$. The geometric realization of $\Sub_p(\square[n])$ is $[0,p]^n$. We obtain easily a well-defined cocubical set $\Sub_p(\square[*])$, and then a colimit-preserving functor 
 \[
 \Sub_p(K) = \int^{[n]\in \square} K_n.\Sub_p(\square[n])\]
with the convention $\Sub_p(\square[0])=\square[0]$. For $p=1$, this functor is the identity. The precubical set $\Sub_3(\square[2])$ is depicted in Figure~\ref{fig:sub3-2}. For the cubical complexes in the sense of \cite{zbMATH01959268} defined in \cite[Section~8.2]{dubut_PhD}, the cubical subdivision described in \cite[Section~8.3.4.2]{dubut_PhD} corresponds to $\Sub_2$. The functor $\Sub_2$ is also mentioned in \cite[Section~2]{MR2521708} with the terminology of \textit{barycentric subdivision}. For any precubical set $K$, there is an isomorphism of directed spaces $\vec{|K|}\iso \vec{|\Sub_p(K)|}$ for all $p\geq 1$ since all involved functors are colimit-preserving and since this holds for $K=\square[n]$ for all $n\geq 0$. 

\bd \label{def:cubical-sbd}
Let $K,L$ be two precubical sets. A \textit{cubical subdivision} $f:K\dashrightarrow L$ consists of an isomorphism of directed spaces $f:\vec{|K|} \iso \vec{|L|}$ such that for each cube $c$ of $K$, the compact $|c|_{geom}$ of the Hausdorff space $|K|_{geom}\iso |L|_{geom}$ is homeomorphic to $|\Sub_{n_c}(c)|_{geom}$ for some $n_c\geq 1$ with $\Sub_{n_c}(c)\subset L$. In particular, $f(K_0)\subset L_0$ which means that $K_0$ can be identified with a subset of $L_0$.
\ed

A cubical subdivision $f:K\dasharrow L$ is not a map of precubical sets unless it is trivial. For example, there exists a cubical subdivision $\square[2]\dasharrow \Sub_3(\square(2])$ (see Figure~\ref{fig:sub3-2}). Besides, it does not give rise to a map of flows $||K||\to ||L||$ in general for an r-realization functor $||-||$. Indeed, the execution path of $\P_{(0,0),(1,1)}||\square[2]||$ depicted in Figure~\ref{fig:dpath-sub} cannot be an execution path of $\P_{(0,0),(1,1)}||\Sub_3(\square[2])||$ because, any r-realization functor being colimit-preserving, the execution path should be the composite of execution paths of some of the nine squares.

\begin{figure}
	\def\n{5}
	\begin{tikzpicture}[black,scale=4,pn/.style={circle,inner sep=0pt,minimum width=4pt,fill=dark-red}]
		\fill [color=gray!15] (0,0) -- (0,1) -- (1,1) -- (1,0) -- cycle;
		\draw[dark-red][-] (0,0) -- (0,1/3) -- (1/3,1/3) -- (1/3,0) -- cycle;
		\draw[dark-red][-] (1/3,0) -- (1/3,1/3) -- (2/3,1/3) -- (2/3,0) -- cycle;
		\draw[dark-red][-] (2/3,0) -- (2/3,1/3) -- (1,1/3) -- (1,0) -- cycle;
		\draw[dark-red][-] (0,1/3) -- (0,2/3) -- (1/3,2/3) -- (1/3,1/3) -- cycle;
		\draw[dark-red][-] (1/3,1/3) -- (1/3,2/3) -- (2/3,2/3) -- (2/3,1/3) -- cycle;
		\draw[dark-red][-] (2/3,1/3) -- (2/3,2/3) -- (1,2/3) -- (1,1/3) -- cycle;
		\draw[dark-red][-] (0,2/3) -- (0,1) -- (1/3,1) -- (1/3,2/3) -- cycle;
		\draw[dark-red][-] (1/3,2/3) -- (1/3,1) -- (2/3,1) -- (2/3,2/3) -- cycle;
		\draw[dark-red][-] (2/3,2/3) -- (2/3,1) -- (1,1) -- (1,2/3) -- cycle;
		\draw[->, line width=0.4mm, color=dark-red, smooth] plot coordinates {(0,0) (1/6,2.5/6) (1,1)};
	\end{tikzpicture}
	\caption{An execution path of $\P_{(0,0),(1,1)}||\square[2]||\backslash \P_{(0,0),(1,1)}||\Sub_3(\square[2])||$}
	\label{fig:dpath-sub}
\end{figure}

\bth \label{thm:sub}
Let $K,L$ be two precubical sets. Consider a cubical subdivision $f:K\dashrightarrow L$. Let $||-||:\square^{op}\set\to \dtop$ be an r-realization functor with $r\in \{q,m,h\}$. For all $\alpha\in K_0$. there is a homotopy equivalence $\P^-_\alpha ||K|| \simeq \P^-_\alpha ||L||$. For all $\alpha\in L_0\backslash K_0$, the space $\P^-_\alpha||L||$ is contractible.
\eth

\bpf
Let $\alpha\in K_0$. Let $h\geq 1$. Consider the set of cubes 
\[
G^h_{\alpha}(K) = \bigg\{c\in \mathcal{C}^-_\alpha(K)\mid n_c=h\bigg\}.
\]
If $c\in G^h_{\alpha}(K)$, then all nonzero dimensional faces of $c$ belong to $G^h_{\alpha}(K)$. Thus
\[
\overline{G}^h_{\alpha}(K) = G^h_{\alpha}(K) \cup \mathcal{V}(G^h_{\alpha}(K))
\]
is a precubical set for all $h\geq 1$ where $\mathcal{V}(G^h_{\alpha}(K))$ is the set of vertices of all cubes of $G^h_{\alpha}(K)$. Moreover, $h\neq h'$ implies $G^h_\alpha(K)\cap G^{h'}_\alpha(K)=\varnothing$. This implies that for all $h\neq h'$, one has $\{\alpha\}\subset \overline{G}^h_{\alpha}(K) \cap \overline{G}^{h'}_{\alpha}(K)\subset K_0$ and $\{\alpha\}\subset\Sub_h(\overline{G}^h_{\alpha}(K)) \cap \Sub_{h'}(\overline{G}^{h'}_{\alpha}(K))\subset K_0$. By definition of $G^h_{\alpha}(K)$, the isomorphism of directed spaces $f:\vec{|K|} \iso \vec{|L|}$ induces an isomorphism of directed spaces $\vec{|\overline{G}^h_{\alpha}(K)|} \iso \vec{|\Sub_h(\overline{G}^h_{\alpha}(K))|}$ for all $h\geq 1$. We then have 
\begin{align*}
	\P^-_\alpha ||K|| & \iso \coprod_{h\geq 1} \P^-_\alpha ||\overline{G}^h_{\alpha}(K)|| \\
	& \simeq \coprod_{h\geq 1} \mathcal{P}^-_\alpha\bigg(\overline{G}^h_{\alpha}(K),\frac{1}{2h}\bigg) \\
	& = \coprod_{h\geq 1} \mathcal{P}^-_\alpha\bigg(\Sub_h(\overline{G}^h_{\alpha}(K)),\frac{1}{2}\bigg)\\
	& \simeq  \coprod_{h\geq 1} \P^-_\alpha ||\Sub_h(\overline{G}^h_{\alpha}(K))||\\
	&\iso \P^-_\alpha ||L||,
\end{align*}
the first isomorphism since $\{\alpha\}\subset \overline{G}^h_{\alpha}(K) \cap \overline{G}^{h'}_{\alpha}(K)\subset K_0$ for $h\neq h'$, the homotopy equivalences by Corollary~\ref{cor:top-flow-homotopy-equivalent}, the equality since $\Sub_h$ increases the natural length by a factor $h$ for all $h\geq 1$, and finally the second isomorphism since $\{\alpha\}\subset\Sub_h(\overline{G}^h_{\alpha}(K)) \cap \Sub_{h'}(\overline{G}^{h'}_{\alpha}(K))\subset K_0$ for $h\neq h'$.

Let $\alpha\in L_0\backslash K_0$. Since $\{\beta\}\subset \overline{G}^h_{\beta}(K) \cap \overline{G}^{h'}_{\beta}(K)\subset K_0$ for all $\beta\in K_0$ and all $h\neq h'$, $\alpha$ belongs to some $\Sub_h(\overline{G}^h_{\beta}(K))$. Thus, one can suppose that $L=\Sub_h(K)$ with $h\geq 2$. If $K=K_{\leq 0}=K_0$, then $K_0=L_0=L$ and the assertion is trivial. If $K=K_{\leq 1}$, then all elements of $L_0\backslash K_0$ are between elements of $K_0$. Thus in this case, $\P^-_\alpha||L||$ is a singleton for all $\alpha\in L_0\backslash K_0$. For all precubical sets $K$, the map $K_{\leq 1}\to K$ is a transfinite composition of pushouts of maps of the form $\de\square[n]\subset \square[n]$ with $n\geq 2$. Thus the map $\Sub_h(K_{\leq 1})\to \Sub_h(K)$ is a transfinite composition of pushouts of maps of the form $\Sub_h(\de\square[n])\subset \Sub_h(\square[n])$ with $n\geq 2$, the functor $\Sub_h$ being colimit-preserving. Let $M:\lambda\to \square^{op}\set$ be a colimit-preserving tower of precubical sets for some ordinal $\lambda$ with $M(0)=K_{\leq 1}$ and $M(\lambda)=\liminj M=K$. Assume that we have a pushout diagram of precubical sets for $0\leq \nu<\lambda$
\[
\begin{tikzcd}[column sep=3em,row sep=3em]
	\arrow[d]\de\square[n]\arrow[r] & M(\nu) \arrow[d] \\
	\square[n]\arrow[r] & \cocartesian M(\nu+1)
\end{tikzcd}
\]
with $n\geq 2$. The functor $\Sub_h$ being colimit-preserving, we obtain the pushout diagram of precubical sets 
\[
\begin{tikzcd}[column sep=3em,row sep=3em]
	\arrow[d]\Sub_h(\de\square[n])\arrow[r] & \Sub_h(M(\nu)) \arrow[d] \\
	\Sub_h(\square[n])\arrow[r] & \cocartesian \Sub_h(M(\nu+1))
\end{tikzcd}
\]
Using Theorem~\ref{thm:G-colimit-preserving}, we obtain the pushout diagram of m-cofibrant spaces 
\[
\begin{tikzcd}[column sep=3em,row sep=3em]
	\mathcal{P}^-(\Sub_h(\de\square[n]),\epsilon)\arrow[d,rightarrowtail]\arrow[r] & \mathcal{P}^-(\Sub_h(M(\nu)),\epsilon) \arrow[d,rightarrowtail] \\
	\mathcal{P}^-(\Sub_h(\square[n]),\epsilon)\arrow[r] & \cocartesian \mathcal{P}^-(\Sub_h(M(\nu+1)),\epsilon)
\end{tikzcd}
\]
where the vertical arrows are m-cofibrations by Corollary~\ref{cor:cubical-cofibration-3}. Either $\alpha\in L_0\backslash K_0$ is not a vertex of the image of the horizontal maps and there is the isomorphism $\mathcal{P}^-_\alpha(\Sub_h(M(\nu)),\epsilon)\iso \mathcal{P}^-_\alpha(\Sub_h(M(\nu+1)),\epsilon)$, or there is the pushout diagram of m-cofibrant spaces 
\[
\begin{tikzcd}[column sep=3em,row sep=3em]
	\mathcal{P}^-_\alpha(\Sub_h(\de\square[n]),\epsilon)\arrow[d,rightarrowtail]\arrow[r] & \mathcal{P}^-_\alpha(\Sub_h(M(\nu)),\epsilon) \arrow[d,rightarrowtail] \\
	\mathcal{P}^-_\alpha(\Sub_h(\square[n]),\epsilon)\arrow[r] & \cocartesian \mathcal{P}^-_\alpha(\Sub_h(M(\nu+1)),\epsilon)
\end{tikzcd}
\]
There are the homeomorphisms \[\mathcal{P}^-_\alpha(\Sub_h(\de\square[n]),\epsilon)\iso N_{n-1}(\epsilon)\hbox{ and }\mathcal{P}^-_\alpha(\Sub_h(\square[n]),\epsilon) \iso N_n(\epsilon).\] Using Proposition~\ref{prop:compacity}, we deduce that the m-cofibrant spaces $\mathcal{P}^-_\alpha(\Sub_h(\de\square[n]),\epsilon)$ and $\mathcal{P}^-_\alpha(\Sub_h(\square[n]),\epsilon)$ are contractible since $n\geq 2$. Consequently, the left vertical map, as well as the right one, are trivial m-cofibrations between contractible m-cofibrant spaces. We then deduce that the map $\mathcal{P}^-_\alpha(\Sub_h(K_{\leq 1}),\epsilon)\to \mathcal{P}^-_\alpha(\Sub_h(K),\epsilon)$ is a transfinite composition of trivial m-cofibrations between contractible m-cofibrant spaces, all involved functors being colimit-preserving. We deduce that $\mathcal{P}^-_\alpha(\Sub_h(K),\epsilon)$ is contractible. The proof is complete since there is the homeomorphism $\P^-_\alpha||\Sub_h(K)||\iso \mathcal{P}^-_\alpha(\Sub_h(K),\epsilon)$ by Corollary~\ref{cor:top-flow-homotopy-equivalent}.
\epf

For completeness, we recall the definition of the branching homology of flows. 

\bd\label{def:hombrdef-flow}  \cite[Definition~6.1]{exbranching}
Let $X$ be a flow. The \textit{branching homology} groups $H_{*}^-(X)$ are defined as follows: 
\begin{enumerate}
	\item for $n\geq 1$, $H_{n+1}^-(X):=H_n(\hop^-X)$
	\item  $H_1^-(X):=\ker(\epsilon)/\im(\partial)$
	\item $H_0^-(X):=\mathbb{Z}[X^0]/\im(\epsilon)$
\end{enumerate}
with the augmentation $\epsilon:\mathbb{Z}[\hop^-X]\to \mathbb{Z}[X^0]$ defined by $\epsilon(\gamma)=\gamma(0)$ and with the map \[\partial:\mathbb{Z}[\top([0,1],\hop^-X)]\longrightarrow \mathbb{Z}[\hop^-X]\] defined by $\partial(f)=f(0)-f(1)$, where $\ker$ is the kernel and $\im$ the image. For any flow $X$, $H_0^-(X)$ is the free abelian group generated by the final states of $X$. The branching homology of $X$ is clearly $X^0$-graded. Let 
\[
H^-_*(X)= \coprod_{\alpha\in X^0} G^\alpha H^-_*(X).
\]
\ed

By \cite[Corollary~6.5]{exbranching}, for every weak equivalence of flows $f:X\to Y$, the morphism of abelian groups $H_n^-(f):H_n^-(X)\to H_n^-(Y)$ is an isomorphism for all $n\geq 0$.

\begin{cor} \label{cor:enfinfin}
	Let $K,L$ be two precubical sets. Let $f:K\dashrightarrow L$ be a cubical subdivision. Then for any $n\geq 0$, there is an isomorphism $H_n^-(||K||)\iso H_n^-(||L||)$ for any r-realization functor $||-||:\square^{op}\set\to \dtop$.
\end{cor}

\bpf
By Corollary~\ref{cor:hop-not-necessary}, there is a homotopy equivalence $\hop^-||K||\simeq \P^-||K||$ for any precubical set $K$. Let $f:K\dashrightarrow L$ be a cubical subdivision. By Theorem~\ref{thm:sub}, for any $\alpha\in K_0$, $f$ induces a homotopy equivalence $\hop^-_\alpha ||K||\simeq \hop^-_{\alpha}||L||$. For any $\alpha\in L_0\backslash K_0$, the branching space $\hop^-_\alpha ||L||$ is contractible by Theorem~\ref{thm:sub}. We deduce that for $n\geq 1$, there are the isomorphisms 
\begin{multline*}
	H_{n+1}^-(||K||) = H_n(\hop^- ||K||) \iso \coprod_{\alpha\in K_0} H_n(\hop^-_\alpha ||K||) \\\iso \coprod_{\alpha\in L_0} H_n(\hop^-_\alpha ||L||) \iso H_n(\hop^- ||L||) = H_{n+1}^-(||L||).
\end{multline*}
None of the states of $\alpha\in L_0\backslash K_0$ is final. Thus all $\alpha \in L_0\backslash K_0$ belongs to the image of the augmentation $\epsilon$. This implies the isomorphism $H_0^-(||K||)\iso H_0^-(||L||)$. For $\alpha\in L_0\backslash K_0$, we have seen that $\hop^-_\alpha||L||$ is contractible, hence path-connected. Thus the quotient of the kernel of the map $\epsilon:\mathbb{Z}[\hop^-_\alpha||L||]\to \mathbb{Z}[\{\alpha\}]$ by the image of the map $\partial:\mathbb{Z}[\top([0,1],\hop^-_\alpha||L||)]\to \mathbb{Z}[\hop^-_\alpha||L||]$ is the trivial abelian group. Hence we deduce that 
\[
H_1^-(||K||)=\coprod_{\alpha\in K_0} G^\alpha H_1^-(||K||) \iso \coprod_{\alpha\in K_0} G^{\alpha} H_1^-(||L||) \iso \coprod_{\alpha\in L_0} G^\alpha H_1^-(|L||) = H_1^-(||L||).
\]
\epf

%\bibliographystyle{../plainurlwithoutprefixDOI} 
%\bibliography{../Bibliotheque}

\end{document}